\documentclass[twoside,reqno]{amsart}

\usepackage{latexsym}

\newcommand{\qdn}{\hspace*{-1.5mm}}
\newcommand{\qqdn}{\hspace*{-2.5mm}}
\newcommand{\xqdn}{\hspace*{-5.0mm}}
\newcommand{\xxqdn}{\hspace*{-10mm}}

\newcommand{\binm}{\binom}
\newcommand{\sbnm}[2]{\Bigl(\!\ba{c}\!#1\!\\#2\ea\!\Bigr)}

\newcommand{\nnm}{\nonumber}
\newcommand{\be}{\begin{equation}}
\newcommand{\ee}{\end{equation}}
\newcommand{\ba}{\begin{array}}
\newcommand{\ea}{\end{array}}
\newcommand{\bmn}{\begin{eqnarray}}
\newcommand{\emn}{\end{eqnarray}}
\newcommand{\bnm}{\begin{eqnarray*}}
\newcommand{\enm}{\end{eqnarray*}}
\newcommand{\bln}{\begin{subequations}}
\newcommand{\eln}{\end{subequations}}

\newtheorem{thm}{Theorem}

\newtheorem{corl}[thm]{Corollary}

\newtheorem{entry}{Entry}

\newcommand{\bbtm}[4]{\bibitem{kn:#1}{#2,}~{#3,}~{#4.}}
\newcommand{\cito}[1]{\cite{kn:#1}}
\newcommand{\citu}[2]{\cite[#2]{kn:#1}}

\begin{document} 
{
\title{Telescoping method, derivative operators and \\harmonic number identities}
\author{$^A$Chuanan Wei, $^B$Dianxuan Gong}


\dedicatory{$^A$Department of Information Technology\\
 Hainan Medical College, Haikou 571101, China\\
 $^B$College of Sciences\\
 Hebei Polytechnic University, Tangshan 063009, China}

\thanks{\emph{Email addresses}: weichuanan@yahoo.com.cn (C. Wei),
gongdianxuan@yahoo.com.cn (D. Gong)}


\begin{abstract}
In terms of the telescoping method, a simple binomial sum is given.
By applying the derivative operators to the equation just mentioned,
we establish several general harmonic number identities including
some known results.
\end{abstract}

\maketitle\thispagestyle{empty}
\markboth{C. Wei,  D. Gong}
         {Telescoping method, derivative operators and harmonic number identities}

\section{Introduction}

For $n\in \mathbb{N}_0$, define the
 harmonic numbers by
\[H_{0}=0\quad \text{and}\quad H_{n}
=\sum_{k=1}^n\frac{1}{k}\quad \text{when}\quad n=1,2,\cdots.\]
 There exist many elegant identities involving harmonic numbers. They
can be found in the papers \cito{andrews}-\cito{zheng}.

For a function $f(x,y)$, define respectively the derivative
operators $\mathcal{D}_x$ and $\mathcal{D}^2_{xy}$ by
 \bnm
&&\xxqdn\mathcal{D}_xf(x,y)=\frac{\partial}{\partial x}f(x,y)\bigg|_{x=0},\\[1mm]
&&\xxqdn\mathcal{D}_{xy}^2f(x,y)=\frac{\partial}{\partial y}
\bigg\{\frac{\partial}{\partial x}f(x,y)\bigg\}\bigg|_{x=y=0}.
 \enm

 Then it is not difficult to show the following two
 derivatives:
 \bnm
&&\mathcal{D}_x\:\binm{s+x}{t}=\binm{s}{t}\big(H_s-H_{s-t}\big),\\
&&\mathcal{D}_{xy}^2\frac{\binm{s+x}{t}}{\binm{u+y}{v}}
=\frac{\binm{s}{t}}{\binm{u}{v}}\big(H_s-H_{s-t}\big)\big(H_{u-v}-H_u\big),
 \enm
where $s,t,u,v\in N_0$ with $t\leq s$ and $v\leq u$.

For a complex sequence $\{\tau_k\}_{k\in \mathbb{Z}}$, define
 the difference operator by
\[\nabla\tau_k=\tau_k-\tau_{k-1}.\]
Then we have the following relation:
 \bnm
\qquad\nabla\frac{\binm{x+k+1}{k}}{\binm{y+k}{k}}=\frac{\binm{x+k}{k}}{\binm{y+k}{k}}\frac{x-y+1}{x+1}.
 \enm
Combining the last equation and the telescoping method:
 \bnm
\sum_{k=1}^{n}\nabla\tau_k=\tau_n-\tau_{0},
 \enm
we obtain the simple binomial sum:
 \bmn \label{source}
\sum_{k=1}^n\frac{\binm{x+k}{k}}{\binm{y+k}{k}}=\frac{\binm{x+n+1}{n}}{\binm{y+n}{n}}\frac{x+1}{x-y+1}-\frac{x+1}{x-y+1}.
 \emn

 By applying the derivative operators to \eqref{source}, several general harmonic
number identities including some known results will be established
in the next two sections.

\section{Harmonic number identities}

\subsection{}

Performing the replacement $y\to y+p$ for \eqref{source} with $p\in
N_0$, we have
  \bnm
\sum_{k=1}^n\frac{\binm{x+k}{k}}{\binm{y+p+k}{k}}=\frac{\binm{x+n+1}{n}}{\binm{y+p+n}{n}}\frac{x+1}{x-y-p+1}-\frac{x+1}{x-y-p+1}.
 \enm
Applying the derivative operator $\mathcal{D}_y$ to the last
equation, we establish the theorem.

\begin{thm}\label{thm-a}
For $x\in \mathbb{C}$ and $p\in \mathbb{N}_0$, there holds the
harmonic number identity:
 \bnm
\quad\sum_{k=1}^n\frac{\binm{x+k}{k}}{\binm{p+k}{k}}H_{p+k}
=\frac{x+1}{x-p+1}\bigg\{\frac{\binm{x+n+1}{n}}{\binm{p+n}{n}}\bigg(H_{p+n}-\frac{1}{x-p+1}\bigg)-H_{p}+\frac{1}{x-p+1}\bigg\}.
 \enm
\end{thm}

Letting $x=p$ in Theorem \ref{thm-a}, we achieve the following
equation.

\begin{corl}\label{corl-a}
For $p\in \mathbb{N}_0$, there holds the harmonic number identity:
 \bnm
\sum_{k=1}^n H_{p+k}=(p+n+1)H_{p+n}-(p+1)H_p-n.
 \enm
\end{corl}

When $p=0$, Corollary \ref{corl-a} reduces to the known result (cf.
\citu{chen}{Equation (2.1)}):
 \bnm
 \sum_{k=1}^nH_k=(n+1)H_{n}-n.
 \enm
Setting $p=n$ and $p=2n$ in Corollary \ref{corl-a} respectively, we
attain the two identities:
 \bnm
 &&\xxqdn\xqdn\,\sum_{k=1}^n H_{n+k}=(2n+1)H_{2n}-(n+1)H_n-n,\\
 &&\xxqdn\xqdn\,\sum_{k=1}^n H_{2n+k}=(3n+1)H_{3n}-(2n+1)H_{2n}-n.
 \enm

Making $x=p+1$ in Theorem \ref{thm-a} and considering the relation:
\[\xxqdn\,\sum_{k=1}^n\frac{p+1+k}{p+1}H_{p+k}=\sum_{k=1}^n\bigg(1+\frac{k}{p+1}\bigg)H_{p+k},\]
we get the following equation by using Corollary \ref{corl-a}.

\begin{corl}\label{corl-b}
For $p\in \mathbb{N}_0$, there holds the harmonic number identity:
 \bnm
\sum_{k=1}^nkH_{p+k}
=\frac{(n-p)(p+n+1)}{2}H_{p+n}+\frac{p(p+1)}{2}H_p-\frac{n(n-2p-1)}{4}.
 \enm
\end{corl}

When $p=0$, Corollary \ref{corl-b} reduces to the result due to Chen
et al. \citu{chen}{Equation (2.2)}:
 \bnm
\xxqdn\qdn \sum_{k=1}^nkH_k=\frac{n(n+1)}{2}H_n-\frac{(n-1)n}{4}.
 \enm
Taking $p=n$ and $p=2n$ in Corollary \ref{corl-b} respectively, we
gain the two identities:
  \bnm
&&\sum_{k=1}^n
kH_{n+k}=\frac{n(n+1)}{2}\Big(H_{n}+\frac{1}{2}\Big),\\
&&\sum_{k=1}^n kH_{2n+k}=\frac{n(3n+1)}{4}(1-2H_{3n})+n(2n+1)H_{2n}.
 \enm

Letting $x=p+2$ in Theorem \ref{thm-a} and considering the relation:
\[\sum_{k=1}^n\frac{(p+1+k)(p+2+k)}{(p+1)(p+2)}H_{p+k}
=\sum_{k=1}^n\bigg\{1+\frac{(2p+3)k}{(p+1)(p+2)}+\frac{k^2}{(p+1)(p+2)}\bigg\}H_{p+k},
\] we achieve the following equation by using Corollaries
\ref{corl-a} and \ref{corl-b}.

\begin{corl}\label{corl-c}
For $p\in \mathbb{N}_0$, there holds the harmonic number identity:
 \bnm
\sum_{k=1}^nk^2H_{p+k}
&&\xqdn\!=\frac{(p+n+1)(2n^2+n-2pn+p+2p^2)}{6}H_{p+n}\\
&&\xqdn\!-\:\frac{p(p+1)(2p+1)}{6}H_p-\frac{n(4n^2-3n-6pn+12p+12p^2-1)}{36}.
 \enm
\end{corl}

When $p=0$, Corollary \ref{corl-c} reduces to the result due to Chen
et al. \citu{chen}{Equation (2.3)}:
 \bnm
\xxqdn\xqdn\qdn\sum_{k=1}^nk^2H_k=\frac{n(n+1)(2n+1)}{6}H_n-\frac{(n-1)n(4n+1)}{36}.
 \enm
When $p=n$, Corollary \ref{corl-c} reduces to the known result (cf.
\citu{chen}{Equation (2.7)}):
 \bnm
\,\quad\sum_{k=1}^nk^2H_{n+k}=\frac{n(n+1)(2n+1)}{6}(2H_{2n}-H_n)-\frac{n(n+1)(10n-1)}{36}.
 \enm
We remark that Chyzak \cito{chyzak} and Schneider \cito{schneider}
proved the last equation by an extension of Zeilberger's algorithm
and Karr's algorithm respectively. Setting $p=2n$ in Corollary
\ref{corl-c}, we attain the identity:
 \bnm
\xxqdn\qdn\sum_{k=1}^nk^2H_{2n+k}&&\xqdn\!=\frac{n(2n+1)(3n+1)}{2}H_{3n}\\
&&\xqdn\!-\:\frac{n(2n+1)(4n+1)}{3}H_{2n}-\frac{n(40n^2+21n-1)}{36}.
 \enm

Choosing $x=-n-1$ in Theorem \ref{thm-a}, we recover the result due
to Sofo \citu{sofo-a}{Corollary 1}:
 \bmn \label{sofo-a}
\xxqdn\xxqdn\xxqdn\:\sum_{k=1}^n(-1)^k\binm{n}{k}\frac{H_{p+k}}{\binm{p+k}{k}}=\frac{-n}{p+n}\bigg(H_{p}+\frac{1}{p+n}\bigg).
 \emn
Choosing $x=-n$ in Theorem \ref{thm-a} and considering the relation:
\[\xxqdn\xqdn\qqdn\sum_{k=1}^n(-1)^k\binm{n-1}{k}\frac{H_{p+k}}{\binm{p+k}{k}}
=\sum_{k=1}^n(-1)^k\binm{n}{k}\bigg(1-\frac{k}{n}\bigg)\frac{H_{p+k}}{\binm{p+k}{k}},\]
we recover, by using \eqref{sofo-a}, the result due to Sofo
\citu{sofo-a}{Corollary 3}:
 \bmn \label{sofo-b}
\quad\:\:\sum_{k=1}^n(-1)^k\binm{n}{k}\frac{kH_{p+k}}{\binm{p+k}{k}}=\frac{n(n^2-n-p^2)}{(p+n)^2(p+n-1)^2}
-\frac{pnH_{p}}{(p+n)(p+n-1)}.
 \emn
Choosing $x=1-n$ in Theorem \ref{thm-a} and considering the
relation:
\[\sum_{k=1}^n(-1)^k\binm{n-2}{k}\frac{H_{p+k}}{\binm{p+k}{k}}
=\sum_{k=1}^n(-1)^k\binm{n}{k}\bigg\{1-\frac{(2n-1)k}{n(n-1)}+\frac{k^2}{n(n-1)}\bigg\}\frac{H_{p+k}}{\binm{p+k}{k}},\]
we recover, by using \eqref{sofo-a} and \eqref{sofo-b}, the result
due to Chu\citu{chu-a}{Corollary 1.8}:
 \bmn
\qqdn\qdn\sum_{k=1}^n(-1)^k\binm{n}{k}\frac{k^2H_{p+k}}{\binm{p+k}{k}}
&&\xqdn\!=\frac{pn(n-p)H_{p}}{(p+n)(p+n-1)(p+n-2)} \nnm\\
\label{chu}&&\xqdn-\:\frac{n^3}{(p+n)^2}+\frac{n(n-1)(2n-1)}{(p+n-1)^2}-\frac{n(n-1)(n-2)}{(p+n-2)^2}.
 \emn
The case $p=0$ of \eqref{sofo-a},  \eqref{sofo-b} and \eqref{chu}
read, respectively, as
 \bmn
&&\sum_{k=1}^n(-1)^k\binm{n}{k}H_{k}=-\frac{1}{n}\quad
\text{where}\quad n>0,\label{sofo-a-spe}\\\label{sofo-b-spe}
&&\sum_{k=1}^n(-1)^k\binm{n}{k}kH_{k}=\frac{1}{n-1}\quad
\text{where}\quad n>1,\\\label{chu-spe}
&&\sum_{k=1}^n(-1)^k\binm{n}{k}k^2H_{k}=-\frac{n}{(n-1)(n-2)}\quad
\text{where}\quad n>2.
 \emn
The last three identities are very beautiful. Subsequently, we shall
display several equations which include also \eqref{sofo-a-spe},
\eqref{sofo-b-spe} and \eqref{chu-spe}.

Making $x\to-x-n-1$, $p\to0$ for Theorem \ref{thm-a}, we get the
following equation.

\begin{corl}\label{corl-d}
For $x\in \mathbb{C}$, there holds the harmonic number identity:
 \bnm
\sum_{k=1}^n(-1)^k\binm{x+n}{k}H_{k}
=(-1)^n\binm{x+n-1}{n}\bigg(H_{n}+\frac{1}{x+n}\bigg)-\frac{1}{x+n}.
 \enm
\end{corl}

When $x=0$, Corollary \ref{corl-d} reduces to \eqref{sofo-a-spe}
exactly. Taking $x=n$ in Corollary \ref{corl-d}, we gain the
identity:
  \bnm
\qquad\sum_{k=1}^n(-1)^k\binm{2n}{k}H_{k}
=(-1)^n\binm{2n-1}{n}\bigg(H_{n}+\frac{1}{2n}\bigg)-\frac{1}{2n}
\quad\text{where}\quad n>0.
 \enm
When $x=-\frac{1}{2}-n$, Corollary \ref{corl-d} reduces to the
result due to Chen et al. \citu{chen}{Example 2.3}:
 \bnm
 \xqdn\sum_{k=1}^n\binm{2k}{k}\frac{H_{k}}{4^k}
=2-\frac{n+1}{4^n}\binm{2n+2}{n+1}+\frac{2n+1}{4^n}\binm{2n}{n}H_n.
 \enm
Letting $x\to x-1$ for Corollary \ref{corl-d} and considering the
relation:
\[\sum_{k=1}^n(-1)^k\binm{x-1+n}{k}H_{k}=\sum_{k=1}^n(-1)^k\bigg(1-\frac{k}{x+n}\bigg)\binm{x+n}{k}H_{k},\]
 we achieve the following equation by using Corollary \ref{corl-d}.

\begin{corl}\label{corl-e}
For $x\in \mathbb{C}$, there holds the harmonic number identity:
 \bnm
\xqdn\sum_{k=1}^n(-1)^k\binm{x+n}{k}kH_{k}
&&\xqdn\!=\frac{(-1)^nnx}{x+n-1}\binm{x+n}{n}\\&&\xqdn\!\times\:
\bigg\{H_{n}+\frac{n^2+(n-1)(x-1)}{n(x+n)(x+n-1)}\bigg\}
+\frac{1}{x+n-1}.
 \enm
\end{corl}

When $x=0$, Corollary \ref{corl-e} reduces to \eqref{sofo-b-spe}
exactly. Setting $x=n$ and $x=-\frac{1}{2}-n$ in Corollary
\ref{corl-e} respectively, we attain the two identities:
  \bnm
\:\quad\sum_{k=1}^n(-1)^k\binm{2n}{k}kH_{k}
=\frac{(-1)^nn^2}{2n-1}\binm{2n}{n}\bigg\{H_{n}+\frac{2n^2-2n+1}{2n^2(2n-1)}\bigg\}+\frac{1}{2n-1},
 \enm
 \bnm
 \xqdn\xxqdn\xxqdn\sum_{k=1}^n\binm{2k}{k}\frac{kH_{k}}{4^k}
=\frac{2n+1}{9}\frac{\binm{2n}{n}}{4^n}(3nH_{n}-2n+6)-\frac{2}{3}.
 \enm

Making $x\to x-1$ for Corollary \ref{corl-e} and considering the
relation:
\[\qdn\xxqdn\sum_{k=1}^n(-1)^k\binm{x-1+n}{k}kH_{k}=\sum_{k=1}^n(-1)^k\bigg(k-\frac{k^2}{x+n}\bigg)\binm{x+n}{k}H_{k},\]
 we get the following equation by using Corollary \ref{corl-e}.

\begin{corl}\label{corl-f}
For $x\in \mathbb{C}$, there holds the harmonic number identity:
 \bnm
&&\xqdn\qqdn\sum_{k=1}^n(-1)^k\binm{x+n}{k}k^2H_{k}
=(-1)^n\frac{x(n^2-n+nx-1)}{(x+n-1)(x+n-2)}\binm{x+n}{n}\\&&\xqdn\qqdn\:\times\:
\bigg\{nH_{n}+\frac{n(2x+2n-3)}{(x+n-1)(x+n-2)}-\frac{n^2+2n-1}{n^2-n+nx-1}\bigg\}
-\frac{x+n}{(x+n-1)(x+n-2)}.
 \enm
\end{corl}

When $x=0$, Corollary \ref{corl-f} reduces to \eqref{chu-spe}
exactly. Taking $x=n$ and $x=-\frac{1}{2}-n$ in Corollary
\ref{corl-f} respectively, we gain the two identities:
  \bnm
&&\xqdn\sum_{k=1}^n(-1)^k\binm{2n}{k}k^2H_{k}
=(-1)^n\frac{n^2(2n+1)}{2(2n-1)}\binm{2n}{n}
\bigg\{H_{n}+\frac{4n^3-8n^2+5n-2}{2n(n-1)(4n^2-1)}\bigg\}
\\&&\qquad\qquad\qquad\qquad\!-\:\frac{n}{(n-1)(2n-1)}
\quad\:\text{where}\quad n>1,\\
 &&\xqdn\sum_{k=1}^n\binm{2k}{k}\frac{k^2H_{k}}{4^k}
=\frac{(2n+1)(3n+2)}{15}\frac{\binm{2n}{n}}{4^n}\bigg(nH_{n}+\frac{8n-6}{9n+6}-\frac{2n}{5}\bigg)+\frac{2}{15}.
 \enm
\subsection{}
Employing the substitutions $x\to y+p$, $y\to x$ for \eqref{source}
with $p\in N_0$, we have
  \bnm
\xqdn\sum_{k=1}^n\frac{\binm{y+p+k}{k}}{\binm{x+k}{k}}=\frac{\binm{y+p+n+1}{n}}{\binm{x+n}{n}}\frac{y+p+1}{y-x+p+1}-\frac{y+p+1}{y-x+p+1}.
 \enm
Applying the derivative operator $\mathcal{D}_y$ to the last
equation, we found the theorem.

\begin{thm}\label{thm-b}
For $y\in \mathbb{C}$ and $p\in \mathbb{N}_0$, there holds the
harmonic number identity:
 \bnm
\xqdn\!\!\sum_{k=1}^n\frac{\binm{p+k}{k}}{\binm{x+k}{k}}H_{p+k}
&&\xqdn\!=\frac{p+1}{p-x+1}\frac{\binm{p+n+1}{n}}{\binm{x+n}{n}}\bigg(H_{p+n+1}
-\frac{1}{p-x+1}\bigg)\\&&\xqdn\!-\:\frac{p+1}{p-x+1}H_{p}+\frac{x}{(p-x+1)^2}.
 \enm
\end{thm}

Of course, Corollaries \ref{corl-a}-\ref{corl-c} can also be implied
by the theorem. Now, we shall derive other several results, which
correspond to Corollaries \ref{corl-d}-\ref{corl-f}, from Theorem
\ref{thm-b}.

Letting $x\to-x-n-1$, $p\to0$ for Theorem \ref{thm-b}, we achieve
 the following equation.

\begin{corl}\label{corl-g}
For $x\in \mathbb{C}$, there holds the harmonic number identity:
 \bnm
\quad\sum_{k=1}^n\frac{(-1)^k}{\binm{x+n}{k}}H_{k}
=\frac{n+1}{x+n+2}\frac{(-1)^n}{\binm{x+n}{n}}\bigg(H_{n+1}-\frac{1}{x+n+2}\bigg)-\frac{x+n+1}{(x+n+2)^2}.
 \enm
\end{corl}

When $x=0$, Corollary \ref{corl-g} reduces to the interesting
identity:
  \bnm
\qquad\:\sum_{k=1}^n\frac{(-1)^k}{\binm{n}{k}}H_{k}
=(-1)^n\frac{n+1}{n+2}\bigg(H_{n+1}-\frac{1}{n+2}\bigg)-\frac{n+1}{(n+2)^2}.
 \enm
Setting $x=n$ and $x=-\frac{1}{2}-n$ in Corollary \ref{corl-g}
respectively, we attain the two identities:
 \bnm
&&\xxqdn\sum_{k=1}^n\frac{(-1)^k}{\binm{2n}{k}}H_{k}
=\frac{(-1)^n}{4\binm{2n}{n}}\bigg(2H_{n+1}-\frac{1}{n+1}\bigg)-\frac{2n+1}{4(n+1)^2},\\
&&\xxqdn\sum_{k=1}^n\frac{4^k}{\binm{2k}{k}}H_{k}
=\frac{2(n+1)}{3}\frac{4^n}{\binm{2n}{n}}\bigg(H_{n+1}-\frac{2}{3}\bigg)-\frac{2}{9}.
 \enm

Making $x\to x+1$ for Corollary \ref{corl-g} and considering the
relation:
\[\sum_{k=1}^n\frac{(-1)^k}{\binm{x+n+1}{k}}H_{k}
=\sum_{k=1}^n\frac{(-1)^k}{\binm{x+n}{k}}\bigg(1-\frac{k}{x+n+1}\bigg)H_{k},\]
 we get the following equation by using Corollary \ref{corl-g}.

\begin{corl}\label{corl-h}
For $x\in \mathbb{C}$, there holds the harmonic number identity:
 \bnm
\sum_{k=1}^n\frac{(-1)^k}{\binm{x+n}{k}}\,kH_{k}
&&\xqdn\!=\bigg\{\frac{1+x+n(x+n+3)}{(x+n+2)(x+n+3)}H_{n+1}-\frac{x+n+1}{(x+n+2)^2}+\frac{x+1}{(x+n+3)^2}\bigg\}
\\&&\xqdn\!\times\:\,\frac{(-1)^n}{\binm{x+n}{n}}\,(n+1)-
\bigg\{\frac{x+n+1}{(x+n+2)^2}-\frac{x+n+2}{(x+n+3)^2}\bigg\}(x+n+1).
 \enm
\end{corl}

When $x=0$, Corollary \ref{corl-h} reduces to the interesting
identity:
  \bnm
\xxqdn\qdn\sum_{k=1}^n\frac{(-1)^k}{\binm{n}{k}}\,kH_{k}
&&\xqdn\!=\bigg\{H_{n+1}+\frac{n}{n^2+3n+1}-\frac{2n+5+(-1)^n}{(n+2)(n+3)}\bigg\}\\&&\xqdn\times\:
\frac{(n+1)(n^2+3n+1)(-1)^n}{(n+2)(n+3)}.
 \enm
Taking $x=n$ and $x=-\frac{1}{2}-n$ in Corollary \ref{corl-h}
respectively, we gain the two identities:
  \bnm
\:\sum_{k=1}^n\frac{(-1)^k}{\binm{2n}{k}}\,kH_{k}
&&\xqdn\!=\bigg\{H_{n+1}+\frac{n}{2n^2+4n+1}-\frac{4n+5}{(2n+2)(2n+3)}\bigg\}
\\&&\xqdn\times\:
\frac{(-1)^n}{\binm{2n}{n}}\frac{2n^2+4n+1}{4n+6}-\frac{(2n+1)(4n^2+6n+1)}{(2n+2)^2(2n+3)^2},\\
&&\xxqdn\xxqdn\xxqdn\:
 \sum_{k=1}^n\frac{4^k}{\binm{2k}{k}}\,kH_{k}
=\frac{2(n+1)(3n+1)}{15}\frac{4^n}{\binm{2n}{n}}\bigg(H_{n+1}-\frac{2}{9n+3}-\frac{2}{5}\bigg)+\frac{2}{225}.
 \enm

Letting $x\to x+1$ for Corollary \ref{corl-h} and considering the
relation:
\[\xxqdn\sum_{k=1}^n\frac{(-1)^k}{\binm{x+n+1}{k}}\,kH_{k}=\sum_{k=1}^n\frac{(-1)^k}{\binm{x+n}{k}}\bigg(k-\frac{k^2}{x+n+1}\bigg)H_{k},\]
 we achieve the following equation by using Corollary \ref{corl-h}.

\begin{corl}\label{corl-i}
For $y\in \mathbb{C}$, there holds the harmonic number identity:
 \bnm
\:\:\:\sum_{k=1}^n\frac{(-1)^k}{\binm{x+n}{k}}\,k^2H_{k}
=(n+1)\frac{(-1)^n}{\binm{x+n}{n}}\big\{\alpha_nH_{n+1}+\beta_n\big\}+\gamma_n.
 \enm
 \bnm
\xxqdn\xxqdn\xxqdn\text{where}\qquad\qquad\qquad\qquad\qquad\quad
&&\xxqdn\xxqdn\xxqdn\xqdn
\alpha_n=\tfrac{1}{x+n+2}-\tfrac{3n+6}{x+n+3}+\tfrac{n^2+5n+6}{x+n+4},\\
&&\xxqdn\xxqdn\xxqdn\xqdn
\beta_n=\tfrac{2x+2n+3}{(x+n+2)^2}-\tfrac{7x+15+2n(x+n+5)}{(x+n+3)^2}+\tfrac{5x+14+n(2x+n+8)}{(x+n+4)^2},\\
&&\xxqdn\xqdn\xxqdn\xxqdn
\gamma_n=\Big\{\tfrac{2x+2n+3}{(x+n+2)^2}-\tfrac{5x+5n+12}{(x+n+3)^2}+\tfrac{3x+3n+10}{(x+n+4)^2}\Big\}(x+n+1).
 \enm
\end{corl}

When $x=0$, Corollary \ref{corl-h} reduces to the interesting
identity:
  \bnm
\xqdn\quad\sum_{k=1}^n\frac{(-1)^k}{\binm{n}{k}}\,k^2H_{k}
&&\xqdn\!=\Big\{\tfrac{n(n^3+7n^2+14n+7)}{(n+2)(n+3)(n+4)}H_{n+1}-
\tfrac{n^6+14n^5+77n^4+208n^3+279n^2+160n+24}{(n+2)^2(n+3)^2(n+4)^2}\Big\}
\\&&\xqdn\!\times\:\:(n+1)(-1)^n-
\tfrac{(n+1)(n^4+5n^3-n^2-28n-24)}{(n+2)^2(n+3)^2(n+4)^2}.
 \enm
Setting $x=n$ and $x=-\frac{1}{2}-n$ in Corollary \ref{corl-h}
respectively, we attain the two identities:
  \bnm
\sum_{k=1}^n\frac{(-1)^k}{\binm{2n}{k}}\,k^2H_{k}
&&\xqdn\!=\Big\{\tfrac{n(n+2)(2n+1)}{2(n+1)(2n+3)}H_{n+1}-
\tfrac{4n^5+28n^4+69n^3+71n^2+28n+3}{4(n+1)^2(n+2)(2n+3)^2}\Big\}
\\&&\xqdn\!\times\:\:\frac{(-1)^n}{\binm{2n}{n}}(n+1)-
\tfrac{n(n+2)(8n^3+8n^2-8n-13)-6}{4(n+1)^2(n+2)^2(2n+3)^2},\\
&&\xxqdn\xxqdn\xxqdn\qdn
 \sum_{k=1}^n\frac{4^k}{\binm{2k}{k}}\,k^2H_{k}
=\Big(\tfrac{15n^2+12n-1}{210}H_{n+1}-\tfrac{225n^2+432n+34}{11025}\Big)
\frac{4^{n+1}}{\binm{2n}{n}}(n+1)+\tfrac{346}{11025}.
 \enm
\section{Further harmonic number identities}
\subsection{}

Performing the replacements $x\to x+p$, $y\to y+q$ for
\eqref{source} with $p,q\in N_0$, we have
  \bnm
\sum_{k=1}^n\frac{\binm{x+p+k}{k}}{\binm{y+q+k}{k}}=\frac{\binm{x+p+n+1}{n}}{\binm{y+q+n}{n}}
\frac{x+p+1}{x-y+p-q+1}-\frac{x+p+1}{x-y+p-q+1}.
 \enm
Applying the derivative operator $\mathcal{D}^2_{xy}$ to the last
equation and using Theorems \ref{thm-a} and \ref{thm-b}, we
establish the theorem.

\begin{thm}\label{thm-c}
For $p,q\in \mathbb{N}_0$, there holds the harmonic number identity:
 \bnm
\sum_{k=1}^n\frac{(p+k)!}{(q+k)!}H_{p+k}H_{q+k}
&&\xqdn\!=\bigg\{(p-q+1)H_{p+n+1}H_{q+n}\!-\!H_{p+n+1}\!-\!H_{q+n}+\frac{2}{p-q+1}\bigg\}\\
&&\xqdn\times\:\frac{(p+n+1)!}{(q+n)!(p-q+1)^2}-\frac{(p+1)!}{q!(p-q+1)^2}
\\&&\xqdn\times\:\bigg\{(p-q+1)H_{p+1}H_{q}-H_{p+1}-H_{q}+\frac{2}{p-q+1}\bigg\}.
 \enm
\end{thm}

Making $q=p$ in Theorem \ref{thm-c} , we get the following equation.

\begin{corl}\label{corl-j}
For $p\in \mathbb{N}_{0}$, there holds the harmonic number identity:
 \bnm
\quad\sum_{k=1}^nH_{p+k}^2=(p+n+1)H_{p+n}^2-(2p+2n+1)H_{p+n}
-(p+1)H_{p}^2+(2p+1)H_{p}+2n.
 \enm
\end{corl}

When $p=0$, Corollary \ref{corl-j} reduces to the known result (cf.
\citu{chen}{Equation (2.8)}):
 \bnm
 \sum_{k=1}^nH_k^2=(n+1)H_{n}^2-(2n+1)H_n+2n.
 \enm
Taking $p=n$ and $p=2n$ in Corollary \ref{corl-j} respectively, we
gain the two identities:
 \bnm
 &&\xqdn\qdn\sum_{k=1}^nH_{n+k}^2=(2n+1)H_{2n}^2-(4n+1)H_{2n}
 -(n+1)H_{n}^2+(2n+1)H_{n}+2n,\\
 &&\xqdn\qdn\sum_{k=1}^nH_{2n+k}^2=(3n+1)H_{3n}^2-(6n+1)H_{3n}
 -(2n+1)H_{2n}^2+(4n+1)H_{2n}+2n.
 \enm
Letting $q=p-1$ in Theorem \ref{thm-c} and considering the relation:
\[\sum_{k=1}^n(p+k)H_{p+k}H_{p+k-1}=\sum_{k=1}^n(p+k)H_{p+k}^2
-\sum_{k=1}^nH_{p+k},\]
 we achieve the following equation by using Corollaries \ref{corl-a} and \ref{corl-j}.

\begin{corl}\label{corl-k}
For $p\in \mathbb{N}_{0}$, there holds the harmonic number identity:
\bnm
 \xqdn\sum_{k=1}^nkH_{p+k}^2&&\xqdn\!=\frac{(n-p)(n+p+1)}{2}H_{p+n}^2
-\frac{n^2-n-1-p(2n+3p+3)}{2}H_{p+n}\\
&&\xqdn+\:\frac{p(p+1)}{2}H_{p}^2-\frac{3p^2+3p+1}{2}H_{p}+\frac{n(n-6p-3)}{4}.
 \enm
\end{corl}

When $p=0$, Corollary \ref{corl-k} reduces to the interesting
identity:
 \bnm
 \quad\sum_{k=1}^nkH_k^2=\frac{n(n+1)}{2}H_{n}^2-\frac{n^2-n-1}{2}H_{n}+\frac{n(n-3)}{4}.
 \enm
Setting $p=n$ and $p=2n$ in Corollary \ref{corl-k} respectively, we
attain the two identities:
 \bnm
 &&\xxqdn\sum_{k=1}^nkH_{n+k}^2=\frac{(2n+1)^2}{2}H_{2n}
 +\frac{n(n+1)}{2}H_{n}^2-\frac{3n^2+3n+1}{2}H_{n}-\frac{n(5n+3)}{4},\\
 &&\xxqdn\sum_{k=1}^nkH_{2n+k}^2=-\frac{n(3n+1)}{2}H_{3n}^2+\frac{15n^2+7n+1}{2}H_{3n}
 +n(2n+1)H_{2n}^2\\&&\qquad\qquad-\:\frac{12n^2+6n+1}{2}H_{2n}-\frac{n(11n+3)}{4}.
 \enm

Making $q=p-2$ in Theorem \ref{thm-c} and considering the relation:
\bnm
\xqdn \sum_{k=1}^n(p+k)(p+k-1)H_{p+k}H_{p+k-2}&&\xqdn\!=\sum_{k=1}^n\{p(p-1)+(2p-1)k+k^2\}H_{p+k}^2\\
&&\xqdn\!-\:\sum_{k=1}^n\{(2p-1)+2k\}H_{p+k},
  \enm
we get the following equation by using Corollaries \ref{corl-a},
\ref{corl-b}, \ref{corl-j} and \ref{corl-k}.

\begin{corl}\label{corl-l}
For $p\in \mathbb{N}_{0}$, there holds the harmonic number identity:
\bnm
 \xqdn\sum_{k=1}^nk^2H_{p+k}^2&&\xqdn\!=\tfrac{(p+n+1)(2n^2+n-2pn+p+2p^2)}{6}H_{p+n}^2-\tfrac{p(p+1)(2p+1)}{6}H_{p}^2
\\&&\xqdn-\tfrac{4n^3-3n^2-6pn^2-n+12pn+12p^2n+3+17p+33p^2+22p^3}{18}H_{p+n}\\[1mm]
&&\xqdn+\:\tfrac{(2p+1)(11p^2+11p+3)}{18}H_{p}+\tfrac{n(8n^2-15n-30pn+25+132p+132p^2)}{108}.
 \enm
\end{corl}

When $p=0$, Corollary \ref{corl-l} reduces to the interesting
identity:
 \bnm
 \xqdn\qdn\sum_{k=1}^nk^2H_k^2=\tfrac{n(n+1)(2n+1)}{6}H_{n}^2-\tfrac{4n^3-3n^2-n+3}{18}H_{n}+\tfrac{n(8n^2-15n+25)}{108}.
 \enm
Taking $p=n$ and $p=2n$ in Corollary \ref{corl-l} respectively, we
gain the two identities:
 \bnm
\xqdn\!\sum_{k=1}^nk^2H_{n+k}^2&&\xqdn\!=\tfrac{n(n+1)(2n+1)}{3}H_{2n}^2-\tfrac{32n^3+42n^2+16n+3}{18}H_{2n}
 -\tfrac{n(n+1)(2n+1)}{6}H_{n}^2\\&&\xqdn+\:\tfrac{(2n+1)(11n^2+11n+3)}{18}H_{n}+\tfrac{n(110n^2+117n+25)}{108},
 \enm
  \bnm
\sum_{k=1}^nk^2H_{2n+k}^2&&\xqdn\!=\tfrac{n(2n+1)(3n+1)}{2}H_{3n}^2-\tfrac{72n^3+51n^2+11n+1}{6}H_{3n}
 -\tfrac{n(2n+1)(4n+1)}{3}H_{2n}^2\\&&\xqdn+\:\tfrac{(4n+1)(44n^2+22n+3)}{18}H_{2n}+\tfrac{n(476n^2+249n+25)}{108}.
 \enm

\subsection{}
Employing the substitutions $x\to -x-p-n-1$, $y\to y+q$ for
\eqref{source} with $p,q\in N_0$, we have
  \bnm
\qquad\sum_{k=1}^n(-1)^k\frac{\binm{x+p+n}{k}}{\binm{y+q+k}{k}}=(-1)^n\frac{\binm{x+p+n-1}{n}}{\binm{y+q+n}{n}}
\frac{x+p+n}{x+y+p+q+n}-\frac{x+p+n}{x+y+p+q+n}.
 \enm
Applying the derivative operator $\mathcal{D}^2_{xy}$ to the last
equation and using Theorem \ref{thm-a}, we found the theorem.

\begin{thm}\label{thm-d}
For $p,q\in \mathbb{N}_0$, there holds the harmonic number identity:
 \bnm
&&\xxqdn\sum_{k=1}^n(-1)^k\binm{p+q+n}{q+k}H_{p+n-k}H_{q+k}\\
&&\xxqdn\:=\bigg\{H_{p}H_{q+n}+\frac{H_{p}}{p+q+n}-\frac{(q+n)H_{q+n}}{p(p+q+n)}
+\frac{p-q-n}{p(p+q+n)^2}\bigg\}\\
&&\xxqdn\:\times\:\binm{p+q+n}{p}\frac{(-1)^np}{p+q+n}-\binm{p+q+n}{q}\frac{p+n}{p+q+n}\\
&&\xxqdn\:\times\:\bigg\{H_{p+n}H_{q}+\frac{H_{p+n}}{p+q+n}-\frac{qH_{q}}{(p+n)(p+q+n)}
+\frac{p-q+n}{(p+n)(p+q+n)^2}\bigg\}.
 \enm
\end{thm}

Letting $q=0$ in Theorem \ref{thm-d} , we achieve the following
equation.

\begin{corl}\label{corl-m}
For $p\in \mathbb{N}_{0}$, there holds the harmonic number identity:
 \bnm
&&\qdn\qqdn\xqdn\sum_{k=1}^n(-1)^k\binm{p+n}{k}H_{k}H_{p+n-k}=(-1)^n\binm{p+n-1}{n}\\
&&\qdn\qqdn\xqdn\:\times\:\bigg\{H_{p}H_{n}+\frac{H_{p}}{p+n}-\frac{nH_{n}}{p(p+n)}
+\frac{p-n}{p(p+n)^2}\bigg\}-\frac{1}{p+n}\bigg\{H_{p+n}+\frac{1}{p+n}\bigg\}.
 \enm
\end{corl}

When $p=0$, Corollary \ref{corl-m} reduces to the interesting
identity:
 \bnm
\xqdn\qdn\sum_{k=1}^n(-1)^k\binm{n}{k}H_{k}H_{n-k}=\frac{(-1)^{n+1}-1}{n}\bigg(H_{n}+\frac{1}{n}\bigg)
\quad\text{where}\:\: n>0.
 \enm
Setting $p=n$ in Corollary \ref{corl-m}, we attain the identity:
 \bnm
\quad\sum_{k=1}^n(-1)^k\binm{2n}{k}H_{k}H_{2n-k}=(-1)^n\binm{2n-1}{n}H_{n}^2\!-\!\frac{1}{2n}\bigg(H_{2n}\!+\!\frac{1}{2n}\bigg)
\quad\text{where}\:\: n>0.
 \enm
Making $p\to p-1$ for Corollary \ref{corl-m} and considering the
relation:
 \bnm
\qqdn\sum_{k=1}^n(-1)^k\binm{p+n-1}{k}H_{k}H_{p+n-1-k}
&&\xqdn\!=\sum_{k=1}^n(-1)^k\binm{p+n}{k}\bigg(1-\frac{k}{p+n}\bigg)H_{k}H_{p+n-k}\\
&&\xqdn\!-\:\frac{1}{p+n}\sum_{k=1}^n(-1)^k\binm{p+n}{k}H_{k},
 \enm
we get the following equation by using Corollaries \ref{corl-d} and
\ref{corl-m}.

\begin{corl}\label{corl-n}
For $p\in \mathbb{N}_{0}$, there holds the harmonic number identity:
 \bnm
&&\xxqdn\sum_{k=1}^n(-1)^k\sbnm{p+n}{k}kH_{k}H_{p+n-k}
=\tfrac{1}{p+n-1}\Big\{H_{p+n}+\tfrac{1}{p+n-1}\Big\}+\tfrac{(-1)^npn}{p+n-1}\sbnm{p+n}{n}\\
&&\xxqdn\:\times\:\Big\{H_{p}H_{n}+\tfrac{n^2+(p-1)(n-1)}{n(p+n)(p+n-1)}H_{p}-\tfrac{n-1}{p(p+n-1)}H_{n}
-\tfrac{p+n-1}{p(p+n)^2}-\tfrac{n-2pn+p-1}{pn(p+n-1)^2}\Big\}.
 \enm
\end{corl}

When $p=0$, Corollary \ref{corl-n} reduces to the interesting
identity:
 \bnm
\xqdn\xqdn\qqdn\qdn\sum_{k=1}^n(-1)^k\sbnm{n}{k}kH_{k}H_{n-k}
&&\xqdn\!=\tfrac{1}{n-1}\Big(H_{n}+\tfrac{1}{n-1}\Big)
\\&&\xqdn\!+\:(-1)^{n}\tfrac{n}{1-n}\Big\{H_{n}+\tfrac{n^2-n+1}{n^2(n-1)}\Big\}
\quad\text{where}\:\: n>1.
 \enm
Taking $p=n$ in Corollary \ref{corl-n}, we gain the identity:
 \bnm
\:\xqdn\sum_{k=1}^n(-1)^k\sbnm{2n}{k}kH_{k}H_{2n-k}
 &&\xqdn\!=\tfrac{1}{2n-1}\Big(H_{2n}+\tfrac{1}{2n-1}\Big)+
 (-1)^{n}\tfrac{n^2}{2n-1}\sbnm{2n}{n}\\
 &&\xqdn\!\times\:\Big\{H_{n}^2+\tfrac{H_n}{2n^2(2n-1)}+\tfrac{4n^2-2n+1}{4n^3(2n-1)^2}\Big\}
\quad\text{where}\:\: n>0.
 \enm

Letting $p\to p-1$ for Corollary \ref{corl-n} and considering the
relation:
 \bnm
\qqdn\sum_{k=1}^n(-1)^k\binm{p+n-1}{k}kH_{k}H_{p+n-1-k}
&&\xqdn\!=\sum_{k=1}^n(-1)^k\binm{p+n}{k}\bigg(k-\frac{k^2}{p+n}\bigg)H_{k}H_{p+n-k}\\
&&\xqdn\!-\:\frac{1}{p+n}\sum_{k=1}^n(-1)^k\binm{p+n}{k}kH_{k},
 \enm
we achieve the following equation by using Corollaries \ref{corl-e}
and \ref{corl-n}.

\begin{corl}\label{corl-o}
For $p\in \mathbb{N}_{0}$, there holds the harmonic number identity:
 \bnm
\qdn\sum_{k=1}^n(-1)^k\sbnm{p+n}{k}k^2H_{k}H_{p+n-k}&&\xqdn\!
=(-1)^n\sbnm{p+n}{n}\Big\{\theta_nH_{n}H_{p}-\lambda_nH_{n}-\mu_nH_{p}-\nu_n\Big\}\\
&&\xqdn\!-\:\varepsilon_nH_{p+n}-\eta_n,
 \enm
 \bnm
\xxqdn\xxqdn\xxqdn\xqdn\text{where}\:\:\:\qquad\qquad\qquad\;\:
 &&\xxqdn\xxqdn\xqdn \theta_n=\tfrac{pn(n^2+pn-n-1)}{(p+n-1)(p+n-2)},\quad\:
 \lambda_n=\tfrac{n(n-1)}{(p+n-1)^2}+\tfrac{n(n-1)(n-2)}{(p+n-2)^2},\\
 &&\xxqdn\xxqdn\xqdn
 \mu_n=\tfrac{p(n^2+2n-1)}{(p+n-1)(p+n-2)}-\tfrac{pn(2p+2n-3)(n^2+pn-n-1)}{(p+n-1)^2(p+n-2)^2},\\
 &&\xxqdn\xxqdn\xqdn \nu_n=\tfrac{2n^2-4n+1}{(p+n-1)^2}+\tfrac{2n(n-1)}{(p+n-1)^3}
 -\tfrac{3n^2-6n+2}{(p+n-2)^2}+\tfrac{2n(n-1)(n-2)}{(p+n-2)^3},\\
 &&\xxqdn\xxqdn\xqdn \varepsilon_n=\tfrac{p+n}{(p+n-1)(p+n-2)},
 \quad\, \eta_n=\tfrac{(p+n)^2-2}{(p+n-1)^2(p+n-2)^2}.\\
 \enm
\end{corl}

When $p=0$, Corollary \ref{corl-o} reduces to the interesting
identity:
 \bnm
\xxqdn\sum_{k=1}^n(-1)^k\sbnm{n}{k}k^2H_{k}H_{n-k}
&&\xqdn\!=(-1)^n\Big\{\tfrac{n+n^2-n^3}{(n-1)(n-2)}H_{n}-\tfrac{n^4-4n^3+6n^2-4n+2}{(n-1)^2(n-2)^2}\Big\}
\\&&\xqdn\!-\:\tfrac{n}{(n-1)(n-2)}\Big\{H_{n}+\tfrac{n^2-2}{n(n-1)(n-2)}\Big\}
\quad\text{where}\:\: n>2.
 \enm
Setting $p=n$ in Corollary \ref{corl-o}, we attain the identity:
 \bnm
\qdn\sum_{k=1}^n(-1)^k\sbnm{2n}{k}k^2H_{k}H_{2n-k}
 &&\xqdn\!=(-1)^n\sbnm{2n}{n}\Big\{\tfrac{2n^3+n^2}{4n-2}H_{n}^2+\tfrac{n}{(2n-1)^2}H_{n}+\tfrac{4n^2-2n+1}{2(2n-1)^3}\Big\}
\\&&\xqdn\!-\:\tfrac{n}{(n-1)(2n-1)}\Big\{H_{2n}+\tfrac{2n^2-1}{n(2n-1)(2n-2)}\Big\}
\quad\text{where}\:\: n>1.
 \enm

\subsection{}
Performing the replacements $x\to x+q$, $y\to -y-p-n-1$ for
\eqref{source} with $p,q\in N_0$, we have
  \bnm
\:\:\:\sum_{k=1}^n(-1)^k\frac{\binm{x+q+k}{k}}{\binm{y+p+n}{k}}=(-1)^n\frac{\binm{x+q+n+1}{n}}{\binm{y+p+n}{n}}
\frac{x+q+1}{x+y+p+q+n+2}-\frac{x+q+1}{x+y+p+q+n+2}.
 \enm
Applying the derivative operator $\mathcal{D}^2_{xy}$ to the last
equation and using Theorem \ref{thm-b}, we establish the theorem.

\begin{thm}\label{thm-e}
For $p,q\in \mathbb{N}_0$, there holds the harmonic number identity:
 \bnm
&&\xqdn\sum_{k=1}^n\frac{(-1)^k}{\binm{p+q+n}{q+k}}H_{p+n-k}H_{q+k}\\
&&\xqdn\:=\bigg\{H_{p}H_{q+n+1}-\frac{H_{p}}{p+q+n+2}-\frac{H_{q+n+1}}{p+q+n+2}
+\frac{2}{(p+q+n+2)^2}\bigg\}\\
&&\xqdn\:\times\:\frac{(-1)^n}{\binm{p+q+n}{p}}\frac{q+n+1}{p+q+n+2}-\frac{1}{\binm{p+q+n}{q}}\frac{q+1}{p+q+n+2}\\
&&\xqdn\:\times\:\bigg\{H_{p+n}H_{q+1}-\frac{H_{p+n}}{p+q+n+2}-\frac{H_{q+1}}{p+q+n+2}
+\frac{2}{(p+q+n+2)^2}\bigg\}.
 \enm
\end{thm}

Making $q=0$ in Theorem \ref{thm-e} , we get the equation.

\begin{corl}\label{corl-p}
For $p\in \mathbb{N}_0$, there holds the harmonic number identity:
 \bnm
&&\xqdn\sum_{k=1}^n\frac{(-1)^k}{\binm{p+n}{k}}H_{k}H_{p+n-k}=\bigg\{H_{p}H_{n+1}-\frac{H_{p}}{p+n+2}-\frac{H_{n+1}}{p+n+2}
+\frac{2}{(p+n+2)^2}\bigg\}\\
&&\xqdn\:\,\times\:\frac{(-1)^n}{\binm{p+n}{p}}\frac{n+1}{p+n+2}-\frac{p+n+1}{(p+n+2)^2}
\bigg\{H_{p+n}-\frac{p+n}{(p+n+1)(p+n+2)}\bigg\}.
 \enm
\end{corl}

When $p=0$, Corollary \ref{corl-p} reduces to the interesting
identity:
 \bnm
\xqdn\xxqdn\sum_{k=1}^n\frac{(-1)^k}{\binm{n}{k}}H_{k}H_{n-k}
=\frac{n+1}{(n+2)^2}\bigg\{\frac{2}{n+2}-H_{n+1}\bigg\}\{1+(-1)^n\}.
 \enm
Taking $p=n$ in Corollary \ref{corl-p}, we gain the identity:
 \bnm
\:\:\,\sum_{k=1}^n\frac{(-1)^k}{\binm{2n}{k}}H_{k}H_{2n-k}
=\frac{(-1)^n}{\binm{2n}{n}}\frac{H_n^2}{2}-\frac{2n+1}{(2n+2)^2}\bigg\{H_{2n}-\frac{n}{(n+1)(2n+1)}\bigg\}.
 \enm

Letting $p\to p+1$ for Corollary \ref{corl-p} and considering the
relation:
 \bnm
\sum_{k=1}^n\frac{(-1)^k}{\binm{p+n+1}{k}}H_{k}H_{p+n+1-k}
&&\xqdn\!=\sum_{k=1}^n\frac{(-1)^k}{\binm{p+n}{k}}\bigg(1-\frac{k}{p+n+1}\bigg)H_{k}H_{p+n-k}\\
&&\xqdn\!+\:\frac{1}{p+n+1}\sum_{k=1}^n\frac{(-1)^k}{\binm{p+n}{k}}H_{k},
 \enm
we achieve the following equation by using Corollaries \ref{corl-g}
and \ref{corl-p}.

\begin{corl}\label{corl-q}
For $p\in N_0$, there holds the harmonic number identity:
 \bnm
\sum_{k=1}^n\frac{(-1)^k}{\binm{p+n}{k}}\,k
H_{k}H_{p+n-k}&&\xqdn\!=(-1)^n\tfrac{n+1}{\binm{p+n}{n}}\Big\{A_nH_{n+1}H_{p}+B_nH_{n+1}-C_nH_{p}+D_n\Big\}\\
&&\xqdn\!-\:E_nH_{p+n+1}+F_n,
 \enm
 \bnm
\xxqdn\xxqdn\xxqdn\xxqdn\text{where}\qquad\:\:\:
 &&\xxqdn A_n=\tfrac{n^2+3n+pn+p+1}{(p+n+2)(p+n+3)},
\qquad\qquad\quad\: B_n=\tfrac{1}{(p+n+2)^2}-\tfrac{n+2}{(p+n+3)^2},\\
&&\xxqdn
C_n=\tfrac{p+n+1}{(p+n+2)^2}-\tfrac{p+1}{(p+n+3)^2},\qquad\quad
D_n=\tfrac{p+n}{(p+n+2)^3}-\tfrac{p-n-1}{(p+n+3)^3}\\
&&\xxqdn
E_n=\tfrac{(p+n+1)^2}{(p+n+2)^2}-\tfrac{(p+n+1)(p+n+2)}{(p+n+3)^2},\:\:
F_n=\tfrac{(p+n)(p+n+1)}{(p+n+2)^3}-\tfrac{(p+n+1)^2}{(p+n+3)^3}.
 \enm
\end{corl}

When $p=0$, Corollary \ref{corl-q} reduces to the interesting
identity:
 \bnm
\sum_{k=1}^n\frac{(-1)^k}{\binm{n}{k}}kH_{k}H_{n-k}&&\xqdn\!=\Big\{\tfrac{3n+5}{(n+3)^3}-\tfrac{2n+2}{(n+2)^3}\Big\}
\Big\{1+(-1)^n\tfrac{(n^2+5n+7)(n^3+3n^2+2n+2)}{n^4+3n^3-6n^2-24n-14}\Big\}\\
&&\xqdn\!-\:\tfrac{(n+1)(n^2+3n+1)}{(n+2)^2(n+3)^2}\Big\{1+(-1)^n\tfrac{n^3+5n^2+6n-1}{n^2+3n+1}\Big\}H_n.
 \enm
Setting $p=n$ in Corollary \ref{corl-q}, we attain the identity:
\bnm
 \xqdn\qqdn\sum_{k=1}^n\frac{(-1)^k}{\binm{2n}{k}}kH_{k}H_{2n-k}&&\xqdn\!
=\frac{(-1)^n}{\binm{2n}{n}}\Big\{\tfrac{2n^2+4n+1}{4n+6}H_n^2
+\tfrac{H_n}{2(2n+3)^2}+\tfrac{4n^2+10n+7}{4(n+1)(2n+3)^3}\Big\}\\
&&\xqdn\!-\:\tfrac{(2n+1)(4n^2+6n+1)}{4(2n^2+5n+3)^2}H_{2n}-\tfrac{2n+1}{4(n+1)^3}+\tfrac{6n+5}{(2n+3)^3}.
 \enm

Making $p\to p+1$ for Corollary \ref{corl-q} and considering the
relation:
 \bnm
\qqdn\sum_{k=1}^n\frac{(-1)^k}{\binm{p+n+1}{k}}kH_{k}H_{p+n+1-k}
&&\xqdn\!=\sum_{k=1}^n\frac{(-1)^k}{\binm{p+n}{k}}\bigg(k-\frac{k^2}{p+n+1}\bigg)H_{k}H_{p+n-k}\\
&&\xqdn\!+\:\frac{1}{p+n+1}\sum_{k=1}^n\frac{(-1)^k}{\binm{p+n}{k}}kH_{k},
 \enm
we get the following equation by using Corollaries \ref{corl-h} and
\ref{corl-q}.

\begin{corl}\label{corl-r}
For $p\in N_0$, there holds the harmonic number identity:
 \bnm
\sum_{k=1}^n\frac{(-1)^k}{\binm{p+n}{k}}\,k^2
H_{k}H_{p+n-k}&&\xqdn\!=(-1)^n\tfrac{n+1}{\binm{p+n}{n}}\Big\{R_nH_{n+1}H_{p}-S_nH_{n+1}-T_nH_{p}-U_n\Big\}\\
&&\xqdn\!+\:V_nH_{p+n+1}-W_n,
 \enm
 \bnm
\xxqdn\xxqdn\xxqdn\xxqdn\text{where}\qquad\qquad\qquad\:\:
&&\xxqdn\xxqdn\!
  R_n=\tfrac{1}{p+n+2}-\tfrac{3n+6}{p+n+3}+\tfrac{(n+2)(n+3)}{p+n+4},\\
&&\xxqdn\xxqdn
 S_n=\tfrac{1}{(p+n+2)^2}-\tfrac{3n+6}{(p+n+3)^2}+\tfrac{(n+2)(n+3)}{(p+n+4)^2},\\
&&\xxqdn\xxqdn
T_n=\tfrac{(p+n+1)^2}{(p+n+2)^2}-\tfrac{(p+1)(2p+2n+3)}{(p+n+3)^2}+\tfrac{(p+1)(p+2)}{(p+n+4)^2},\\
&&\xxqdn\xxqdn\!
U_n=\tfrac{2p+2n+2}{(p+n+2)^3}-\tfrac{2n^2+7n+2pn+7p+9}{(p+n+3)^3}+\tfrac{3n+2pn+5p+8}{(p+n+4)^3}\\
&&\xxqdn\xxqdn
V_n=\Big\{\tfrac{2p+2n+3}{(p+n+2)^2}-\tfrac{5p+5n+12}{(p+n+3)^2}+\tfrac{3p+3n+10}{(p+n+4)^2}\Big\}(p+n+1),\\
&&\xxqdn\xxqdn\qdn
W_n=\Big\{\tfrac{2p+2n+2}{(p+n+2)^3}-\tfrac{5p+5n+9}{(p+n+3)^3}+\tfrac{3p+3n+8}{(p+n+4)^3}\Big\}(p+n+1).
 \enm
\end{corl}

When $p=0$, Corollary \ref{corl-r} reduces to the interesting
identity:
 \bnm
&&\qdn\xqdn\sum_{k=1}^n\frac{(-1)^k}{\binm{n}{k}}k^2H_{k}H_{n-k}=(n+1)H_{n+1}\\
&&\qdn\xqdn\:\times\:\,
\Big\{\Big(\tfrac{2n+3}{(n+2)^2}\!-\!\tfrac{5n+12}{(n+3)^2}\!+\!\tfrac{3n+10}{(n+4)^2}\Big)
-(-1)^n\Big(\tfrac{1}{(n+2)^2}\!-\!\tfrac{3n+6}{(n+3)^2}+\tfrac{n^2+5n+6}{(n+4)^2}\Big)\Big\}\\
 &&\qdn\xqdn\:-\,\,
 \Big\{\Big(\tfrac{2n+2}{(n+2)^3}\!-\!\tfrac{5n+9}{(n+3)^3}\!+\!\tfrac{3n+8}{(n+4)^3}\Big)
+(-1)^n\Big(\tfrac{2n+2}{(n+2)^3}\!-\!\tfrac{2n^2+7n+9}{(n+3)^3}\!+\!\tfrac{3n+8}{(n+4)^3}\Big)\Big\}(n+1).
 \enm
Taking $p=n$ in Corollary \ref{corl-r}, we gain the identity:
 \bnm
&&\qdn\:\xqdn\sum_{k=1}^n\frac{(-1)^k}{\binm{2n}{k}}k^2H_{k}H_{2n-k}
=\frac{(-1)^nn}{\binm{2n}{n}}\Big\{\tfrac{2n^2+5n+2}{4n+6}H_n^2
+\tfrac{H_n}{(2n+3)^2}+\tfrac{4n^2+10n+7}{(2n+2)(2n+3)^3}\Big\}\\
&&\xqdn\:\,-\:\,\tfrac{(2n+1)(4n^4+10n^3-n^2-14n-6)}{(n+1)^2(2n+3)^2(2n+4)^2}H_{2n}
+\tfrac{n(8n^6+12n^5-132n^4-512n^3-731n^2-459n-103)}{4(n+1)^3(n+2)^3(2n+3)^3}.
 \enm



\end{document}